\documentclass[11pt, oneside, a4paper]{article}
\usepackage[utf8]{inputenc}
\usepackage[T1]{fontenc}
\usepackage{amsmath}
\usepackage{amsthm,enumerate} 
\usepackage{amssymb}
\usepackage{listings}
\usepackage{url}
\usepackage{hyperref}

\newtheorem{thm}{Theorem}[section]
\newtheorem{lemma}[thm]{Lemma}

\newtheorem{theorem}[thm]{Theorem}
\newtheorem{cor}[thm]{Corollary}
\newtheorem{prop}[thm]{Proposition}
\newtheorem{notation}[thm]{Notation}
\newtheorem{remark}[thm]{Remark}
\newtheorem{l+d}[thm]{Lemma und Definition}
\newtheorem{k+d}[thm]{Corollary und Definition}

\newtheorem{example}[thm]{Example}

\newtheorem{bem+de}[thm]{Defintion+Bemerkung}

\newcommand{\erz}[1]{\langle #1\rangle}
\DeclareMathOperator{\PG}{PG}

\title{On the largest independent sets in the Kneser graph on chambers of PG(4,q)}

\author{Philipp Heering\footnote{Justus-Liebig-Universität, Mathematisches Institut, Arndtstraße 2, D-35392 Gießen, Germany, philipp.heering@math.uni-giessen.de } }
\date{April 2024}

\begin{document}

\maketitle

\begin{center}
\subsection*{Abstract}
\end{center}

Let $\Gamma_4$ be the graph whose vertices are the chambers of the finite projective $4$-space $\PG(4,q)$, with two vertices being adjacent if the corresponding chambers are in general position. 
For $q\geq 749 $ we show that $\alpha:=(q^2+q+1)(q^3+2q^2+q+1)(q+1)^2$ is the independence number of $\Gamma_4$ and the geometric structure of independent sets with $\alpha$ vertices is described. \\

\textbf{Keywords:} Erd\H{o}s-Ko-Rado problem, $q$-analog of generalized Kneser graph, 
independence number\\
\textbf{MSC (2020):} 
05B25, 
05C69, 
51E20 

\section{Introduction}

We consider \emph{chambers} $\{P,\ell,\pi, S\}$ of a finite projective $4$-space $\PG(4,q)$ consisting of a point $P$, a line $\ell$, a plane $\pi$ and a solid $S$ that are mutually incident. Two chambers $\{ P_1,\ell_1,\pi_1,S_1\}$ and $\{ P_2,\ell_2,\pi_2,S_2\}$ are called \emph{opposite} if the following conditions are satisfied: $P_1$ does not lie in $S_2$, and $P_2$ does not lie in $S_1$, and $\ell_1$ does not meet $\pi_2$, and $\ell_2$ does not meet $\pi_1$. We denote chambers as tuples $(P,\ell,\pi,S)$. Similarly, we write $(\ell,\pi)$, or $(\ell,\pi,S)$ for subflags of a chamber.\\
 Let $\Gamma_4$ be the graph whose vertices are the chambers of $\PG(4,q)$, with two vertices of $\Gamma_4$ being adjacent if the corresponding chambers are opposite. We call $\Gamma_4$ the \emph{Kneser graph on chambers of $\PG(4,q)$} and study maximal independent sets of $\Gamma_4$. The investigation of independent sets of Kneser graphs is the famous Erd\H{o}s-Ko-Rado problem. Using a suitable notion of intersection these problems have been studied in many different settings, see for example \cite{ OriginalEKR, FRANKL1986228, EKR_integer_sequences, IVANOV2023113363, Thechromaticnumberoftwo} 
and \cite{GodsilMeagher} for an overwiev.  \\
 For projective spaces of odd dimension $d$, the independence number of the Kneser graph on chambers of $\PG(d,q)$ has been determined in \cite{AlgebraicApproach}. For $d=3$ large maximal independent sets have been studied in \cite{maximalindependentsets3d-HeeringMetsch}.
For flags of projective spaces, the analogue problem has also been studied in many different cases, see for example \cite{EKRpointplaneflags4d, cocliquesonlineplane4d, pointhyperplaneflags, werner_metsch}. 
 
We turn our attention back to independent sets of $\Gamma_4$. Utilizing algebraic methods and the Hoffman-ratio bound, it was shown in \cite{AlgebraicApproach} that
 $$\dfrac{(q^4+q^3+q^2+q+1)(q^3+q^2+q+1)(q^2+q+1)(q+1)}{1+q^{(5/2)}}$$
 is an upper bound for the independence number of $\Gamma_4$. However, there are no known examples that meet this bound.

We use a geometric approach. The advantage is that in addition to the independence number of $\Gamma_4$, we also get the structure of the largest maximal independent sets and some information on the second largest maximal independent sets. The downside of our approach is that our proof is only significant if $q$ is large. Our main result is the following.

\begin{theorem} \label{T: main theorem}
		A maximal independent set of the Kneser graph on chambers of $\PG(4,q)$ contains 
 	$$ \alpha:=(q^2+q+1)(q^3+2q^2+q+1)(q+1)^2, $$
 	or at most 
 	$$ \beta:= q^7 + 4 q^6 + 133 q^5 + 290 q^4 + 261 q^3 + 195 q^2 + 75 q + 36$$
 	or at most
 	$$ \gamma:= 750 q^6 + 2379 q^5 + 3914 q^4 + 3910 q^3 + 2755 q^2 + 1225 q + 473  $$
 	chambers. In particular $\alpha$ is the independence number of $\Gamma_4$, for $q\geq 749 $. If $q\geq 749$ and $M$ is a maximal independent set with $|M|>max\{ \beta, \gamma\}$, then $M$ has the structure described in Example \ref{E: maximal sets}.
\end{theorem}

 The structure described in Example \ref{E: maximal sets} is essentially a blowup of maximal independent sets of the Kneser graph one line-plane flags of $\PG(4,q)$ that was studied in \cite{cocliquesonlineplane4d}.\\
Let $M$ be a maximal independent set of $\Gamma_4$. If an incident line-plane pair $(\ell,\pi)$ of $\PG(4,q)$ occurs in a chamber in $M$, we call it an $M$-coflag. 
The idea used to structure the proof of Theorem \ref{T: main theorem} is the following: Let $(\ell,\pi)$ be an $M$-coflag that occurs in $n$ chambers in $M$. We call $n$ the weight of $(\ell,\pi)$. In this case $n$ is one of the following numbers: $1$, $2$, $q+1$, $2q+1$, $(q+1)^2$. This is discussed at length in Section \ref{Section: Weight of a line-plane flag}. The concept of weight that is used here, is similar to the one used in \cite{maximalindependentsets3d-HeeringMetsch}.

\section{Preliminaries} \label{Section: Weight of a line-plane flag}
Consider the finite $d$-dimensional projective space $\PG(d,q)$.
 Let $S_i$ be a subspace of dimension $i$ for all $0\leq i\leq d-1$ with $S_{i}\subseteq S_{i+1}$ for $0\leq i \leq d-2$. Denoting chambers as tuples, we get that $(S_0,\ldots,S_{d-1})$ is a chamber of $\PG(d,q)$. 

 Two chambers $C_1=(S_0,\ldots,S_{d-1})$ and $C_2=(R_0,\ldots,R_{d-1})$ of $PG(d,q)$ are called opposite if  for any two elements $S_i\in C_1$ and $R_j\in C_2$ the subspaces $S_i$ and $R_j$ have either no common point or if they span the entire space $\PG(d,q)$.
 This is equivalent to $S_i\cap R_{d-1-i}=\emptyset$ for all $0\leq i\leq d-1$.
 
By $\Gamma_d$ we denote the \emph{Kneser graph on chambers of $\PG(d,q)$}, that is the vertices of $\Gamma_d$ are the chambers of $\PG(d,q)$ and two vertices are adjacent if the corresponding chambers are opposite. 
 
 \begin{remark} \label{R: independent set of projective plane}
 	Let $M$ be a maximal independent set of $\Gamma_2$. Then one of the two following cases occurs.
 	\begin{itemize}
 		\item[a)] We have $|M|=2q+1$ and there is a point-line flag $(P,\ell)$ in $PG(2,q)$, such that all chambers of $M$ have $\ell$ as their line, or $P$ as their point
		\item[b)] We have $|M|=3$ and the points of the chambers are not collinear.
 	\end{itemize}
  \end{remark}
 
Chambers $(S_0,\ldots,S_{d-1})$ of $\PG(d,q)$ are also denoted as triples $(P,f,H)$, with $(P,f,H)=(S_0,\ldots,S_{d-1})$ if and only if $P=S_0$ and $f=(S_1,\ldots,S_{d-2})$ and $H=S_{d-1}$. In this case we call $f$ a \emph{coflag}. 
Two coflags $f=(S_1,\ldots,S_{d-2})$ and $g=(R_1,\ldots,R_{d-2})$ are called \emph{f-opposite} if  for any two elements $S_i\in f$ and $R_j\in g$ the subspaces $S_i$ and $R_j$ have either no common point or if they span the entire space $\PG(d,q)$. This condition is in some sense inherited from the definition of oppositeness for chambers.
If two coflags $f$  and $g$ are non-f-opposite, we say that $f$ is \emph{non-f-opposite} to $g$.
 In the case $d=4$ coflags are incident line-plane pairs. 

\subsection{Weight of a coflag}

The notation of this subsection is adapted from \cite{maximalindependentsets3d-HeeringMetsch}.

\begin{notation} 
\begin{itemize} Let $M$ be an independent set of the Kneser graph on chambers of $\PG(d,q)$.
	\item[a)] Chambers in $M$ are called $M$-chambers. 
	\item[b)] If $f$ is a coflag, such that $M$ contains a chamber that contains $f$, we call $f$ an $M$-coflag. We also say that $f$ is contained in $M$.
	\item[c)] Let $f$ be an $M$-coflag. The number of $M$-chambers that contain $f$ is called the $M$-weight, or (if there is no risk of confusion) simply weight of $f$.
	\item[d)] Let $(P,f)$ be an incident point-coflag pair. The number of $M$-chambers that contain $P$ and $f$ is called the $M$-weight or simply weight of $(P,f)$. The same terminology is used for coflag-hyperplane pairs.
	\item[e)] Let $P$ be a point and let $f$ be a coflag. We say that $P$ is a point of $f$, or on $f$, if there is a hyperplane $H$, such that $(P,f,H)$ is a chamber. A similar notation is used for cofllags and hyperplanes.
\end{itemize}
\end{notation}

\begin{lemma} \label{L: R contains one point}
	Let $(P,f,H)$ and $(Q,g,R)$ be two chambers of $\PG(d,q)$ with coflags $f$ and $g$ that are f-opposite. Then $R$ contains exactly one of the $q+1$ points on $f$.
\end{lemma}

\begin{proof}
	Let $S_1,\ldots,S_{d-2},S'_1,\ldots,S'_{d-2}$ be subspaces of $\PG(d,q)$, such that
	$f=( S_1,\ldots,S_{d-2} )$ and $g=( S'_1,\ldots,S'_{d-2})$. Since $g$ and $f$ are f-opposite we have $S_1\cap S'_{d-2}=\emptyset$. Since $S'_{d-2}\subseteq R$, we have that $R\cap S_1$ is one of the $q+1$ points on the line $S_1$.
\end{proof}

The following two results are generalizations of Lemma 2.2 and Proposition 2.3 of \cite{maximalindependentsets3d-HeeringMetsch}.

\begin{lemma} \label{L: weight of triple}
Let $M$ be a maximal independent set of chambers of $\PG(d,q)$. Every incident coflag-hyperplane pair $(f,H)$ has weight $0$, $1$, or $q+1$. If a pair $(f,H)$ has weight $q+1$, every $M$-chamber contains a coflag that is non-f-opposite to $f$ or has its point on $H$.
\end{lemma}

\begin{proof}
	Let $(f,H)$ be an incident coflag-hyperplane pair with weight at least two. In this case there are two chambers $(P_1,f,H)$ and $(P_2,f,H)$ in $M$ with $P_1\neq P_2$. Let $P$ be any point of $f$. We show that the chamber $(P,f,H)$ is opposite to no chamber of $M$. Then maximality of $M$ implies that $(P,f,H)\in M$ and since this applies to all points $P$ of $f$, it follows that $(f,H)$ has weight $q+1$.

Let $(Q,g,R)$ be any chamber of $M$. If $g$ is non-f-opposite to $f$, then $(Q,g,R)$ and $(P,f,H)$ are not opposite. Suppose that $g$ and $f$ are f-opposite. We show that $Q\in H$, which implies the second statement. Lemma \ref{L: R contains one point} states that $R$ contains only one of the $q+1$ points on $f$. Without loss of generality, we assume $P_1\notin R$. Since $(P_1,f,H)$ and $(Q,g,R)$ are not opposite, this implies that $Q\in H$. Hence $(Q,g,R)$ and $(P,f,H)$ are not opposite.
\end{proof}

For incident point-coflag pairs we get the dual statement.

\begin{prop}  \label{P: weight of lines}
Let $M$ be a maximal independent set of chambers of $\PG(d,q)$ and let $f$ be a coflag.
Then $f$ has weight $0$, $1$, $2$, $q+1$, $2q+1$ or $(q+1)^2$. Moreover, the following hold.
\begin{enumerate}
\item[a)] The coflag $f$ has weight $(q+1)^2$ if and only if it is non-f-opposite to all $M$-coflags.
\item[b)] If $f$ has weight $2q+1$, then $f$ is incident with a point $P$ and a hyperplane $H$ such that every $M$-chamber $(Q,g,R)$ with a coflag $g$ that is f-opposite to $f$ satisfies $Q\in H$ and $P\in R$. In particular, every $M$-coflag $g$ with the upper mentioned property has weight $1$ and every chamber that contains $f$, also contains $P$ or $H$.
\item[c)] If $f$ has weight $q+1$, then one of the following two cases occurs.
\begin{enumerate}
\item[(i)] There exists a point $P$ on $f$, such that the $M$-chambers of $f$ are the $q+1$ chambers that contain $P$, and $f$. In this case every $M$-chamber $(Q,g,R)$ with a coflag $g$ that is f-opposite to $f$, satisfies $P\in R$.
\item[(ii)] There exists a hyperplane $H$ on $f$, such that the $M$-chambers of $f$ are the $q+1$ chambers that contain $f$ and $H$. In this case every $M$-chamber $(Q,g,R)$ with a coflag $g$ that is f-opposite to $f$, satisfies $Q\in H$.
    \end{enumerate}
\item[d)] If $f$ has weight two and $(P_i,f, H_i)$, $i=1,2$ are the two chambers of $M$ containing $f$, then $P_1\neq P_2$ and $H_1\neq H_2$. In this case, every chamber $(Q,g,R)$ in $M$ with a coflag $g$ that is f-opposite to $f$, satisfies either $P_1\in R$ and $Q\in H_2$, or $P_2\in R$ and $Q\in H_1$.
\end{enumerate}
\end{prop}
\begin{proof}

If no $M$-coflag is f-opposite to  $f$, then none of the $(q+1)^2$ chambers on $f$ are opposite to any $M$-chamber, so maximality of $M$ implies that these $(q+1)^2$ chambers all belong to $M$, so $f$ has weight $(q+1)^2$.

From now on we assume that $M$ contains a chamber $(Q_0,g_0,R_0)$, with $g_0$ f-opposite to $f$. Lemma \ref{L: R contains one point} states that $R_0$ meets the line of $f$ in a point $P_0$, and that the hyperline of $f$ and $Q_0$ span a hyperplane $H_0$. If $(P,f,H)$ is any chamber of $M$ containing $f$, then it is non-opposite to $(Q_0,g_0,R_0)$ and hence $P\in R_0$ or $Q_0\in H$, that is $P=P_0$ or $H=H_0$. Lemma \ref{L: weight of triple} implies that the number of chambers of $M$ containing $P_0$ and $f$ is either $0$, $1$ or $q+1$, and that the number of chamber of $M$ containing $H_0$ and $f$ is either $0$, $1$ or $q+1$. If $q+1$ occurs in any of the two situations mentioned above, then $(P_0,f, H_0)$ lies in $M$. It follows immediately that the number of $M$-chambers on $f$ is $0$, $1$, $2$, $q+1$ or $2q+1$. It remains to prove b)-d).

b) We assume that $f$ lies in exactly $2q+1$ $M$-chambers. Then these are all chambers containing $P_0$ and $f$ or containing $f$ and $H_0$. 
 For the second statement in (b) suppose that $(Q,g,R)$ is an $M$-chamber with $g$ f-opposite to $f$. Since $(Q,g,R)$ is non-opposite to the chambers containing $P_0$ and $f$, we have $P_0\in R$. Dually we get $Q\in H_0$.

c) We assume that $f$ has weight $q+1$. Then these are either the $q+1$ chambers containing $P_0$ and $f$ or the $q+1$ chambers containing $f$ and $H_0$. Suppose first that these are the $q+1$ chambers containing $P_0$ and $f$. Then Lemma \ref{L: weight of triple} proves the second part of $(c)(i)$.
Suppose now that the $M$-chambers that contain $f$ are the $q+1$ chambers containing $f$ and $H_0$. This is dual to the previous situation and therefore we are in situation $(c)(ii)$.

d) Assume that $f$ has weight two and let $(P_i,f, H_i)$, $i=1,2$ be the two chambers of $M$ containing $f$. It is impossible that $P_1=P_2$, since otherwise Lemma \ref{L: weight of triple} implies that $(P_1,f)$ has weight $q+1$ and hence $f$ would have weight at least $q+1$. Dually we have $H_1\not=H_2$. Suppose that $(Q,g,R)$ is a chamber of $M$ with $g$ f-opposite to $f$. Since $(P_i,f, H_i)$ and $(Q,g, R)$ are non-opposite, we have $Q\in H_i$ or $P_i\in R$. Since $g$ and $f$ are f-opposite, Lemma \ref{L: R contains one point} states that $R$ contains at most one of the points $P_1$ and $P_2$. As $H_1\not=H_2$, then $Q$ lies in at most one of these planes. This proves (d).
\end{proof}

\begin{cor}  \label{C: weight 1 properties}
Let $M$ be a maximal independent set of chambers of $\PG(d,q)$.
Let $(P,f,H)$ be in $M$ and suppose that $f$ has weight one.
	Then there are chambers $(Q_1,g_1,R_1)$ and $(Q_2,g_2,R_2)$ in $M$ that satisfy the following properties: $f$ is f-opposite to $g_1$ and $g_2$,  $P\in R_1$, $Q_1\notin H$ and $P\notin R_2$, $Q_2\in H$.
\end{cor}

\begin{proof}
	 There are $M$-coflags that are f-opposite to $f$, otherwise $f$ would have weight $(q+1)^2$. Let us assume that there is no chamber $(Q_1,g_1,R_1)$ in $M$ that satisfies: $f$ is f-opposite to $g_1$, $P\in R_1$, $Q_1\notin H$. This means that all chambers $(Q,g,R)$ in $M$ with a coflag $g$ that is f-opposite to $f$, satisfy $Q\in H$. In this case the pair $(f,H)$ would have weight $q+1$, this implies that $f$ has weight $\geq q+1$, which is a contradiction. Therefore, we can find a chamber $(Q_1,g_1,R_1)$ in $M$ that satisfies the conditions: $f$ is f-opposite to $g_1$, $P\in R_1$, $Q_1\notin H$ and with the same argument we can also find a chamber $(Q_2,g_2,R_2)$ in $M$ that satisfies the conditions: $f$ is f-opposite to $g_2$, $P\notin R_2$, $Q_2\in H$.
\end{proof}

\subsection{Coflags in $\PG(4,q)$}
\label{Section: line-plane flags}

The coflags of $\PG(4,q)$ are line-plane flags $(\ell,\pi)$. 
Let $\Gamma_4(1,2)$ be the graph whose vertices are the coflags of $\PG(4,q)$ with two coflags $(\ell_1,\pi_1)$ and $(\ell_2,\pi_2)$ being adjacent if they are f-opposite, that is $\ell_1\cap \pi_2=l_2\cap \pi_1=\emptyset$. We call $\Gamma_4(1,2)$ the \emph{Kneser graph on line-plane flags of $\PG(4,q)$}. 
We will see that maximal independent sets of $\Gamma_4$ are induced by maximal independent sets of $\Gamma_4(1,2)$. Maximal independent sets of $\Gamma_4(1,2)$ are well studied.

\begin{thm} [\cite{cocliquesonlineplane4d}] \label{T: B and B q^5}
Let $N$ be a largest maximal independent set of $\Gamma_4(1,2)$. Then $|N|=(q^2+q+1)(q^3+2q^2+q+1)$
	 and $N$ is one of the independent sets described in Example \ref{E: line-plane sets max}.
\end{thm}

\begin{example} \label{E: line-plane sets max}
\begin{itemize}
	\item[a)] Let $S_0$ be a solid and $P_0\in S_0$ a point. Let $N$ be the set of all line-plane flags $(\ell,\pi)$ that satisfy $\pi \subseteq S_0$, or $\ell=\pi\cap S_0$ and $P_0\in l$.
	\item[b)] Let $S_0$ be a solid and $\pi_0\subset S_0$ a plane. Let $N$ be the set of all line-plane flags $(\ell,\pi)$ that satisfy $\pi \subseteq S_0$, or $\ell=\pi\cap S_0$ and $\ell\subseteq \pi_0$.
	\item[c)] Let $S_0$ be a solid and $P_0\in S_0$ a point. Let $N$ be the set of all line-plane flags $(\ell,\pi)$ that satisfy $P_0\in\ell$, or $\pi = \erz{ \ell, P_0 }$ and $\pi\subseteq  S_0$.
	\item[d)] Let $\ell_0$ be a line and $P_0\in \ell_0$ a point. Let $N$ be the set of all line-plane flags $(\ell,\pi)$ that satisfy $P_0\in \ell$, or $\pi = \erz{ \ell, P_0 }$ and $\ell_0\subseteq \pi$.
\end{itemize}
\end{example}

Furthermore we have the following Hilton-Milner type theorem.

\begin{thm} [\cite{DHAESELEER2022103474}] \label{T: Hilton-Milner for line-planes}
	Let $N$ be an independent set of $\Gamma_4(1,2)$ that is not a subset of any of the sets described in Example \ref{E: line-plane sets max}. Then $|N|\leq 4q^4+9q^3+4q^2+q+1$
\end{thm}

For later use, we need the following technical lemmata. 

\begin{lemma} \label{L: z 2 line-plane flags opposite}
Let $\lambda:=q^2+q+1$ and
	let $(\ell_1,\pi_1)$ and $(\ell_2,\pi_2)$ be two line-plane flags in $\PG(4,q)$ with $\ell_1\neq \ell_2$ and $\pi_1\neq \pi_2$. The number of line-plane flags $(g,\tau)$ that are non-f-opposite to $(\ell_1,\pi_1)$ and $(\ell_2,\pi_2)$ is at most $z_B:=\lambda (8q^4+14q^3+16q^2+10q+6)$.
\end{lemma}

\begin{proof}
 A line-plane flag $(g,\tau)$ that is non-f-opposite to $(\ell_1,\pi_1)$ and $(\ell_2,\pi_2)$ has to satisfy one of the following conditions: 
	\begin{itemize}
		\item[a)]  $g\cap \pi_1\neq \emptyset$ and $g\cap \pi_2\neq \emptyset$.
		\item[b)] $g\cap \pi_1\neq \emptyset$ and $\ell_2\cap \tau \neq \emptyset$.
		\item[c)] $\ell_1\cap \tau \neq \emptyset$ and $g\cap \pi_2 \neq \emptyset$.
		\item[d)] $\ell_1\cap \tau \neq \emptyset$ and $\ell_2\cap \tau \neq \emptyset$.
	\end{itemize}
	
	We obtain an upper bound for the number of chambers that satisfy each condition, in every case. Condition a) and d) are dual, furthermore note that condition b) and c) are inherently symmetrical. Therefore, the bound we obtain for a), will also be used for d) and the bound we obtain for b), will also be used for c). Every line $g$ is incident with at most $\lambda$ planes $\tau$.

 	\begin{itemize}
 		\item[a)] First, we find an upper bound for the number $n$ of lines $g$ that meet both $\pi_1$ and $\pi_2$. Then clearly $n\cdot \lambda$ is our desired bound. Since $\pi_1\neq \pi_2$ the subspace $\pi_1\cap \pi_2$ can be a  point, or a line.\\
 		 If $\pi_1$ and $\pi_2$ meet in a line, we have that $(q^2\cdot q^2)+(q+1)(q^3+q^2+q)+1$ distinct lines meet both $\pi_1$ and $\pi_2$.\\
 		  If $\pi_1$ and $\pi_2$ meet in a point, we have that $(q^2+q)(q^2+q)+(q^3+q^2+q+1)$ distinct lines meet both $\pi_1$ and $\pi_2$. Since $q\geq 2$, we have $n\leq (q^2\cdot q^2)+(q+1)(q^3+q^2+q)+1$ in any case and therefore $((q^2\cdot q^2)+(q+1)(q^3+q^2+q)+1)\lambda$ is our bound.
 		\item[b)] There are three possible arrangements: Either (i) $\ell_2\subseteq \pi_1$, or (ii) the subspace $\ell_2\cap \pi_1$ is a point $P$, or (iii) $\ell_2\cap \pi_1=\emptyset$.
 		\begin{itemize}
 			\item[(i)] Let $\ell_2\subseteq \pi_1$. If $g\cap \pi_1$ is one of the $q^2$ points not on $\ell_2$ and $\tau\cap \ell_2$ is one of the $q+1$ points on $\ell_2$, then these two points span a line and $\tau$ is one of the $\lambda-1$ planes $\neq \pi_1$ incident with this line. Furthermore $g$ is one of the $q+1$ lines on $\tau$ incident with $g\cap \pi_1$.
 			In every other case the line $g$ meets $\ell_2$, that is $g$ is one of the $1+(q+1)(q^3+q^2+q)$ lines that meets $\ell_2$ and $\tau$ is one of the $\lambda$ planes incident with $g$. Therefore, our bound in this case 
 					is $q^2(q+1)\lambda (q+1) +(1+(q+1)(q^3+q^2+q))\lambda$.
 			\item[(ii)] 
 					Let $\ell_2\cap \pi_1$ be a point $P$. If $g\cap \pi_1$ is one of the $q^2+q$ points not on $\ell_2$ and $\tau\cap \ell_2$ is one of the $q+1$ points on $\ell_2$, then these two points span a line and $\tau$ is one of the $\lambda-1$ planes $\neq \pi_1$ incident with this line. Furthermore $g$ is one of the $q+1$ lines on $\tau$ incident with $g\cap \pi_1$. In every other case the line $g$ meets $\ell_2$, that is $g$ is one of the $1+(q+1)(q^3+q^2+q)$ lines that meets $\ell_2$ and $\tau$ is one of the $\lambda$ planes incident with $g$. Therefore, our bound in this case is $(q^2+q)(q+1)\lambda (q+1) +(1+(q+1)(q^3+q^2+q))\lambda$.
 			\item[(iii)] 			Let $\ell_2\cap \pi_1=\emptyset$. If $g$ does not meet $\ell_2$, then  $g\cap \pi_1$ is one of the $q^2+q+1$ points on $\pi_1$ and $\tau\cap \ell_2$ is one of the $q+1$ points on $\ell_2$. These two points span a line and $\tau$ is one of the $\lambda$ planes incident with this line. Furthermore $g$ is one of the $q+1$ lines on $\tau$ incident with $g\cap \pi_1$. In every other case the line $g$ meets $\ell_2$, that is $g$ is one of the $1+(q+1)(q^3+q^2+q)$ lines that meets $\ell_2$ and $\tau$ is one of the $\lambda$ planes incident with $g$. Therefore, our bound in this case is $(q^2+q+1)(q+1)\lambda (q+1) +(1+(q+1)(q^3+q^2+q))\lambda$.
 		\end{itemize}
 		Clearly $(q^2+q+1)(q+1)\lambda q +(1+(q+1)(q^3+q^2+q))\lambda$ is an upper bound in (i), (ii) and (iii).
 		 	\end{itemize}
In conclusion we obtain an upper bound of
 \begin{align*}
 	\begin{split}
 		 & 2 \cdot \Bigl( ((q^2\cdot q^2)+(q+1)(q^3+q^2+q)+1)\lambda \Bigr) \\
 		+& 2\cdot \Bigl( (q^2+q+1)(q+1)\lambda (q+1) +(1+(q+1)(q^3+q^2+q))\lambda \Bigr) \\
  		=& \lambda (8q^4+14q^3+16q^2+10q+6)
 	\end{split}
 \end{align*}
 as claimed.
\end{proof}

\begin{remark}
	It is clear that the bound obtained in Lemma \ref{L: z 2 line-plane flags opposite} is not optimal. Counting more carefully would allow us to improve the bounds, however we decide against this at this point and also in Subsection \ref{Subsection: weight one or two} for an optimized reading experience.
\end{remark}

\begin{lemma} \label{L: z 2 line-plane flags opposite, Case A}
	Let $(\ell_1,\pi_1)$ and $(\ell_2,\pi_2)$ be two line-plane flags in $\PG(4,q)$ with $\ell_1\neq \ell_2$ and $\pi_1\neq \pi_2$. Furthermore let $P_0$ be a point on $(\pi_1\cap\pi_2)\setminus (\ell_1\cup \ell_2)$. The number of line-plane flags $(g,\tau)$ with $P_0\in \tau\setminus g$ that are non-f-opposite to $(\ell_1,\pi_1)$ and $(\ell_2,\pi_2)$ is at most $z_A:=2q^5+12q^4+11q^3+8q^2+3q+3$.
\end{lemma}

\begin{proof}
The proof works analogue to the proof of Lemma \ref{L: z 2 line-plane flags opposite}.
 A line-plane flag $(g,\tau)$ that is non-f-opposite to $(\ell_1,\pi_1)$ and $(\ell_2,\pi_2)$ has to satisfy one of the following conditions: 
	\begin{itemize}
		\item[a)]  $g\cap \pi_1\neq \emptyset$ and $g\cap \pi_2\neq \emptyset$.
		\item[b)] $g\cap \pi_1\neq \emptyset$ and $\ell_2\cap \tau \neq \emptyset$.
		\item[c)] $\ell_1\cap \tau \neq \emptyset$ and $g\cap \pi_2 \neq \emptyset$.
		\item[d)] $\ell_1\cap \tau \neq \emptyset$ and $\ell_2\cap \tau \neq \emptyset$.
	\end{itemize}
	
	We obtain an upper bound for the number of chambers that satisfy each condition, in every case. For b) and c) we use the same bound. Since $\pi_1\neq \pi_1$ and $P_0\in (\pi_1\cap\pi_2)\setminus (\ell_1\cup \ell_2) $ the line $\ell_1$ cannot be contained in $\pi_2$ and vice versa.
	Note that for every line $g$ the plane $\tau$ is determined by $P_0$. Therefore, we define $\lambda:=1$.
 	\begin{itemize}
 		\item[a)] First, we find an upper bound $n$ for the number of lines $g$ that meet both $\pi_1$ and $\pi_2$. This is analogue to Case a) of Lemma \ref{L: z 2 line-plane flags opposite} and we obtain a bound of $((q^2\cdot q^2)+(q+1)(q^3+q^2+q)+1)\lambda$.
 		\item[b)] There are two possible arrangements: Either (ii) the subspace $\ell_2\cap \pi_1$ is a point $P$, or (iii) $\ell_2\cap \pi_1=\emptyset$.
 		\begin{itemize}
 		
 			\item[(ii)] 
 					Let $\ell_2\cap \pi_1$ be a point $P$. 
 					Let $g\cap \pi_1$ be one of the $q^2+q-1$ points $\neq P,P_0$ on $\pi_1$. Let $Q$ be one of the $q+1$ points on $\ell_2$, such that $g\cap \pi_1, P_0$ and $Q$ are not collinear. Then $\tau$ is determined by these three points and $g$ is one of the $q^2$ lines on $\tau$ that is not incident with $P_0$. If the three points would be collinear, then $\tau$ would be one of the $q^2+q+1$ planes incident with the line $\pi_1\cap \pi_2$.
						In every other case the line $g$ meets $\ell_2$, that is $g$ is one of the $1+(q+1)(q^3+q^2+q)$ lines that meets $\ell_2$ and $\tau$ is one of the $\lambda$ planes incident with $g$. Therefore, our bound in this case is $(q^2+q-1)(q+1)q^2+ (q^2+q+1)q^2 +(1+(q+1)(q^3+q^2+q))\lambda$.
 			\item[(iii)] 			Let $\ell_2\cap \pi_1=\emptyset$. If $g\cap \pi_1$ is one of the $q^2+q$ points $\neq P_0$ on $\pi_1$ and $\tau\cap \ell_2$ is one of the $q+1$ points on $\ell_2$, then these two points span a line that is not incident with $P_0$ and $\tau$ is determined. Furthermore $g$ is one of the $q+1$ lines on $\tau$ incident with $g\cap \pi_1$. In every other case the line $g$ meets $\ell_2$, that is $g$ is one of the $1+(q+1)(q^3+q^2+q)$ lines that meets $\ell_2$ and $\tau$ is one of the $\lambda$ planes incident with $g$. Therefore, our bound in this case is $(q^2+q)(q+1)(q+1) +(1+(q+1)(q^3+q^2+q))\lambda$.
 		\end{itemize}
 		Clearly $ (q^2+q-1)(q+1)q^2+ (q^2+q+1)q^2 +(1+(q+1)(q^3+q^2+q))$ is an upper bound in (ii) and (iii).
 		\item[d)] First, we find an upper bound $n$ for the number of planes $\tau$ that meet both $\ell_1$ and $\ell_2$ for $\ell_1\neq \ell_2$. Then $g$ can be any of the $q^2$ lines on $\tau$ that is not incident with $P_0$, hence $nq^2$ gives a bound.\\
 		 Since $\ell_1\neq \ell_2$, the lines do not meet, or they meet in a point. 
 		If the lines meet in a point, the point $P_0$ is not in the span of $\ell_1$ and $\ell_2$.\\
 		 If $\ell_1\cap \ell_2$ is a point $P$, there are at most $ q^2+q+1$ planes incident with $P$ and $P_0$. Let us assume that $\tau$ is a plane that is incident with $\ell_1$ and $\ell_2$, but not with $P$. Then $\tau$ contains one of the $q^2$ lines of $\erz{\ell_1,\ell_2}$ that is not incident with $P$. Therefore $\tau$ is one of the $\lambda$ planes that is incident with such a line and $P_0$. In conclusion $n\leq (q^2+q+1)+q^2$ in this case.\\
 		 If $\ell_1\cap \ell_2=\emptyset$, there is at most one point $Q$ on $\ell_2$ such that the span of $Q$ and $\ell_1$ contains $P_0$. If $\tau$ contains $QP_0$ then $\tau$ is any of the $(q^2+q+1)$ planes incident with $QP_0$. If $\tau$ does not contain $QP_0$ then 
the plane $\tau$ contains one of the $(q+1)^2$ lines that meet both $\ell_1$ and $\ell_2$. Then $\tau$ is determined by this line and $P_0$. Therefore $n\leq (q^2+q+1)+ (q+1)^2$ in this case.\\
 		 We have $n\leq (q^2+q+1)+ (q+1)^2$ in both cases and therefore $((q^2+q+1)+ (q+1)^2)q^2$ is our bound.
 	\end{itemize}
In conclusion we obtain an upper bound of
 \begin{align*}
 	\begin{split}
 		 & ((q^2\cdot q^2)+(q+1)(q^3+q^2+q)+1) \\
 		+& 2\cdot \Bigl(  (q^2+q-1)(q+1)q^2+ (q^2+q+1)q^2 +(1+(q+1)(q^3+q^2+q)) \Bigr) \\
 		+& ((q^2+q+1)+ (q+1)^2)q^2 \\
  		=& 2q^5+12q^4+11q^3+8q^2+3q+3
 	\end{split}
 \end{align*}
 as claimed.
\end{proof}

\section{Maximal independent sets of $\Gamma_4$}
\label{Section: The largest maximal independent set}

Before we begin with the proof of Theorem \ref{T: main theorem}, we give a class of examples of large maximal independent sets of $\Gamma_4$.

\begin{example} \label{E: maximal sets}
Let $N$ be a largest maximal independent set of the Kneser graph one line-plane flags of $\PG(4,q)$, that is $N$ is one of the sets described in Example \ref{E: line-plane sets max}. Let $M$ be the set of all chambers whose line-plane flag is in $N$. 
Then $M$ is a maximal independent set of the Kneser graph on chambers of $\PG(4,q)$. Theorem \ref{T: B and B q^5} yields $|N|=(q^2+q+1)(q^3+2q^2+q+1)$ and every line-plane flag is contained in at most $(q+1)^2$ chambers (see Proposition \ref{P: weight of lines}). Therefore, we get 
 $|M|=(q^2+q+1)(q^3+2q^2+q+1)(q+1)^2$.
\end{example}

Note that in the above example every coflag of $M$ has weight $(q+1)^2$.\\
 On the other hand, if we take a set $M$ of pairwise non-opposite chambers and assume that every coflag in $M$ has weight $(q+1)^2$, Proposition \ref{P: weight of lines} yields that no two $M$-coflags are f-opposite. Using Theorem \ref{T: B and B q^5} and Theorem \ref{T: Hilton-Milner for line-planes} this implies that either $M$ is as in Example \ref{E: maximal sets} and
\begin{align*}
	|M|=(q^2+q+1)(q^3+2q^2+q+1)(q+1)^2
\end{align*}
or
\begin{align*}
	|M|\leq (4q^4+9q^3+4q^2+q+1)(q+1)^2.
\end{align*}

The remainder of this section is devoted to the case in which a set $M$ of pairwise non-opposite chambers contains a coflag of weight $<(q+1)^2$.

\begin{prop} \label{P: upper bound for (q+1)^2-lines}
Let $M$ be an independent set of $\Gamma_4$ that contains coflags of weight $<(q+1)^2$. Then one of the following two cases occurs.
\begin{itemize}
	\item[A:] The number of $M$-coflags with weight $(q+1)^2$ is at most $(q^3+q^2+q+1)(q^2+q+1)+q^2(q+1)$ and (A1) these include all line-plane flags whose lines contain a fixed point $P_0$, or dually (A2) all line-plane flags whose planes are contained in a fixed solid $S_0$.
	\item[B:] The number of $M$-coflags with weight $(q+1)^2$, is at most $\delta:=4q^4+9q^3+4q^2+q+1$.
\end{itemize}
\end{prop}

\begin{proof}
 Proposition \ref{P: weight of lines} states that no $M$-coflag is f-opposite to any $M$-coflag of weight $(q+1)^2$. 
 Let $N$ be the set of all $M$-coflags with weight $(q+1)^2$.
 If $N$ contains $\delta +1$ elements, the Theorems \ref{T: B and B q^5} and \ref{T: Hilton-Milner for line-planes} imply that $N$ is contained in one of the structures described in Example \ref{E: line-plane sets max}. \\
 Since a) and b) of Example \ref{E: line-plane sets max} are dual to c) and d) we assume without loss of generality that $N$ is contained in the structure of c) or d).
In case c) there is a point $P_0$ and a solid $S_0$ on $P_0$, such that all line-plane flags in $N$ have their line on $P_0$, or their plane in $S_0$. 
In case d) there is a point $P_0$ and a line $\ell_0$ on $P_0$, such that all line-plane flags in $N$ have their line on $P_0$, or their plane on $\ell_0$. \\
Since $|N|>\delta$ the set $N$ contains more than $q^2(q^2+q+1)+(q^2+q+1)^2$ flags. This implies that $N$ contains at least $(q^2+q+1)^2+1$ line-plane flags with lines that contain $P_0$. The lines of these flags cannot all be contained in a solid. This implies that all line-plane flags, that are non-f-opposite to all coflags in $N$, have their plane on $P_0$. Therefore all $(q^3+q^2+q+1)(q^2+q+1)$ line-plane flags with a line on $P_0$ are in $N$.\\
Now assume that $N$ contains $(q^3+q^2+q+1)(q^2+q+1)+q^2(q+1)+1$ coflags, i.e. $q^2(q+1)+1$ coflags with a line that is not on $P_0$. Then these $q^2(q+1)+1$ coflags all have their plane on $P_0$. However, they cannot all have their plane in a fixed solid and on a line. In other words, we are either in case c) of Example \ref{E: line-plane sets max}, or in case d). However, in this case 
 the only line-plane flags that are non-f-opposite to all line-plane flags in this structure, are also in this structure. Hence all $M$-coflags have to be in this structure (and therefore they are non-f-opposite to all $M$-coflags) and this implies that they have weight $(q+1)^2$. Therefore, we get that $M$ does not contain coflags of weight $<(q+1)^2$.
 Since $M$ contains coflags of weight $<(q+1)^2$, we obtain a contradiction.
\end{proof}

From now on we refer to Case A, A1, A2 and B of Proposition \ref{P: upper bound for (q+1)^2-lines} simply as Case A, A1, A2 or B.

\begin{remark} \label{R: raez}
	In Case A1, all $M$-coflags have to have their plane on a fixed point $P_0$ and all coflags with a line on $P_0$ have weight $(q+1)^2$. Therefore, a line of weight $<(q+1)^2$ already determines the line-plane flag. In Case A2 the plane already determines the line. 
\end{remark}
 
Lemma \ref{L: z 2 line-plane flags opposite} concerns the Case B and Lemma \ref{L: z 2 line-plane flags opposite, Case A} concerns the Case A1, dually we also get the same result in Case A2.
Now we obtain an upper bound on the number of $M$-coflags with weight $<(q+1)^2$ in Case A and B.

\subsection{Bound for coflags with weight $q+1$ or $2q+1$}
\label{Subsection: q+1}

\begin{notation}
Let $M$ be an independent set of $\Gamma_4$.
	Let $(\ell,\pi)$ be an $M$-coflag of weight $q+1$, or $2q+1$. If there is a point $Q$, so that $(Q,\ell,\pi)$ has $M$-weight $q+1$, we call $(\ell,\pi)$ a $P_M$-pair. If there is a solid $R$, so that $(\ell,\pi,R)$ has $M$-weight $q+1$, we call $(\ell,\pi)$ an $S_M$-pair.
\end{notation}

Proposition \ref{P: weight of lines} yields that every $M$-coflag of weight $q+1$ is a $P_M$-pair or an $S_M$-pair, furthermore it implies that an $M$-coflag has weight $2q+1$, if and only if it is a $P_M$-pair and an $S_M$-pair.

\begin{lemma} \label{L: S-pairs non-opposite}
For an independent set $M$ of $\Gamma_4$, no two $S_M$-pairs are f-opposite.
\end{lemma}

The proof is adapted from Lemma 2.9 of \cite{maximalindependentsets3d-HeeringMetsch}.

\begin{proof}
	Let $(\ell_1,\pi_1)$ and $(\ell_2,\pi_2)$ be $S_M$-pairs. Let $P_1,\ldots,P_{q+1}$ be the points on $\ell_1$ and let $Q_1,\ldots,Q_{q+1}$ be the points on $\ell_2$. Furthermore let $S_1$ and $S_2$ be the solids, so that the triples $(\ell_1,\pi_1,S_1)$ and $(\ell_2,\pi_2,S_2)$ have weight $q+1$. Then all of the  following chambers are in $M$:
	 $$(P_1,\ell_1,\pi_1,S_1),\ldots,(P_{q+1},\ell_1,\pi_1,S_1), (Q_1,\ell_2,\pi_2,S_2),\ldots,(Q_{q+1},\ell_2,\pi_2,S_2).$$ 
	Not let us assume that $(\ell_1,\pi_1)$ and $(\ell_2,\pi_2)$ are f-opposite. This implies for any two chambers $(P_{i},\ell_1,\pi_1,S_1)$ and $(Q_j,\ell_2,\pi_2,S_2)$ that $P_i\in S_2$, or $Q_j\in S_1$. Without loss of generality, we assume that two distinct points $P_1$, $P_2$ are in $S_2$. This implies that $\ell_1$ is in $S_2$. Therefore $\ell_1$ meets $\pi_2$ in a point. This is a contradiction. 
\end{proof}

For $P_M$-pairs, we have the dual statement: No two $P_M$-pairs are f-opposite.

\begin{prop} \label{P: q+1} 
Let $M$ be an independent set of $\Gamma_4$ and assume that we are in Case B.
	The number of chambers in $M$ with $M$-coflags of weight $q+1$, or $2q+1$ is at most $2\cdot(q^5+3q^4+4q^3+4q^2+2q+1)(q+1)$.
\end{prop}

\begin{proof}
	First, we show that the number of $S_M$-pairs in $M$ is at most $a:=q^5+3q^4+4q^3+4q^2+2q+1$. In Lemma \ref{L: S-pairs non-opposite} we showed that the line-plane flags of two $S_M$-pairs are non-f-opposite. This implies that the number of $S_M$-pairs is bounded by the number of pairwise non-f-opposite line-plane flags. Theorem \ref{T: B and B q^5} yields that the number of pairwise non-opposite line-plane flags is at most $a$. Therefore, the number of $S_M$-pairs in $M$ is at most $a$.\\
	For $P_M$-pairs we have the dual statement: There are at most $a$ $P_M$-pairs in $M$.\\
	Now we need to show that $2a$ is an upper bound for the number of $M$-coflags that are an $S_M$-pair or a $P_M$-pair. If no $S_M$-pair is a $P_M$-pair this is clearly the case. If an $S_M$-pair is also a $P_M$-pair the total number of chambers with $S_M$- and $P_M$-pairs is actually smaller than $2a$, since an $M$-coflag with weight $2q+1$ is a $P_M$-pair and an $S_M$-pair and $2q+1<2(q+1)$.
\end{proof}

\begin{prop} \label{P: q+1 CASE A}
	Let $M$ be an independent set of $\Gamma_4$ and assume that we are in Case A.
	The number of chambers in $M$ with $M$-coflags of weight $q+1$, or $2q+1$ is at most $2\cdot (q^4+q^3+q^2)(q+1)$.
\end{prop}
\begin{proof}
In Case A2 there is a solid $S_0$, such that all $M$-coflags of weight $<(q+1)^2$ have their line in $S_0$ but their plane not in $S_0$.
	We only consider Case A2, since Case A1 is dual.
	 Let $(\ell_1,\pi_1)$ and $(\ell_2,\pi_2)$ be two line-plane flags with $\ell_i \subseteq S_0$ and $\pi_i\not\subseteq S_0$ for a fixed solid $S_0$ and $i=1,2$. If $(\ell_1,\pi_1)$ and $(\ell_2,\pi_2)$ are non-f-opposite, the lines have to meet in a point. The number of lines in $S_0$ that pairwise meet is, according to \cite{FRANKL1986228}, at most $q^2+q+1$. Every line is incident with $q^2$ planes not in $S_0$. Using similar arguments as in Proposition \ref{P: q+1}, we get that the number of chambers with coflags with weight $q+1$, or $2q+1$ is bound by $2\cdot (q^4+q^3+q^2)(q+1)$.
\end{proof}

\subsection{Bound for coflags with weight one or two}
\label{Subsection: weight one or two}

In this subsection we always consider an independent set $M$ of $\Gamma_4$ that contains coflags of weight $<(q+1)^2$.
First, we want to show that the number of $M$-chambers with a coflag of weight one or two, whose point is in a fixed solid, is small. We need the following technical lemmata for the proof.

\begin{lemma} \label{L: points in a solid preparation} 
Let $M$ be an independent set of $\Gamma_4$ and $\lambda:=q^2+q+1$. Let $\sigma$ be a plane and let $Y$ be the set of all chambers $(P,\ell,\pi,S)$ of $M$ for which $(\ell,\pi)$ has $M$-weight at most two and for which $P\in \sigma$. 
In Case B we have $|Y|\leq y_B:=\lambda (8 q^4 + 19 q^3 + 15 q^2 + 5 q + 3)$.
\end{lemma}

\begin{proof} Note that every line in $\PG(4,q)$ is incident with $\lambda$ planes.

Let $Y_3$ be the set of all chambers in $Y$, whose line is contained in $\sigma$. 
Let $Y_2$ be the set of chambers in $Y$, whose plane meets $\sigma$ in a line. Let $Y_1$ be the set of chambers in $Y$, whose plane meets $\sigma$ only in a point. 
Since $Y_1\cup Y_2$ contains all chambers in $Y$ with a line that is not in $\sigma$, we have 
$Y=Y_1\cup Y_2\cup Y_3$. We obtain bounds for $|Y_1|, |Y_2|, |Y_3|$, therefore we assume without loss of generality that $Y_1$, $Y_2$ and $Y_3$ are not empty.

There are $q^2+q+1$ lines in $\sigma$ and each line is incident with $\lambda$ planes. The $Y_3$-coflags have at most weight two, therefore $|Y_3|\leq 2(q^2+q+1)\lambda$.
	
Now we want to find an upper bound for $|Y_1\cup Y_2|$. The point of every chamber of $Y_1\cup Y_2$ is in $\sigma$, but the lines of the chambers of $Y_1\cup Y_2$ are not contained in $\sigma$.
If we assume that a $Y_1\cup Y_2$-coflag $f$ has weight $2$ Proposition \ref{P: weight of lines} yields that there are two distinct points $P_1$ and $P_2$ on $f$ such that $(P_1,f)$ and $(P_2,f)$ have $M$-weight $>0$. But in this case only one of the two points $P_1$, $P_2$ is in $\sigma$, therefore only one of the point-coflag pairs occurs in a chamber of $Y$.
 Therefore, we can count the $Y_1\cup Y_2$-coflags in order to count the chambers in $Y_1\cup Y_2$.

Let $(P,\ell,\pi,S)$ be a chamber of $Y_2$, note that $\pi\neq \sigma$, since otherwise the line $\ell$ would be in $\sigma$. 
There are $q^2+q+1$ possibilities for $\ell\cap \sigma$. Each of these points is incident with $q+1$  lines in $\sigma$. Each of these lines is incident with $\lambda$ planes $\pi$. The line $\ell$ has to be one of the $q+1$ lines in $\pi$ incident with $\ell\cap \sigma$. In conclusion we obtain $|Y_2|\leq (q^2+q+1)(q+1)^2\lambda$.

Now we find an upper bound for $|Y_1|$. 
Let $(P,\ell,\pi,S)$ be a chamber in $Y_1$, i.e. $\pi\cap \sigma=P$. 
Since $Y_1$ is an independent set, all chambers $(Q,g,\tau,R)$ in $Y_1$ have to be non-opposite to $(P,\ell,\pi,S)$.
We find an upper bound for the number of chambers $(Q,g,\tau,R)$ in $Y_1$ that a) contain a line-plane pair that is non-f-opposite to $(\ell,\pi)$, b) have a point $Q\in S$, or c) have a solid $R$ on $P$.
This yields an upper bound for $|Y_1|$. Recall that it is sufficient to count the $Y_1$-coflags.
\begin{itemize}
	\item[a)] We find an upper bound for the number $(\#a)$ of coflags $(g,\tau)$ of chambers in $Y_1$ that are non-f-opposite to $(\ell,\pi)$. 
	
		 The sets $\pi\cap \sigma $ and $ \tau\cap \sigma $ are points, since the corresponding chambers are in $Y_1$.  If $g\cap \sigma=\ell\cap \sigma$, the flags are clearly non-opposite. The point $\ell\cap \sigma$ is incident with $q^3+q^2$ lines $g$ that are not in $\sigma$ and each line is incident with $\lambda$ planes. This yields that there are $(q^3+q^2)\lambda$ flags $(g,\tau)$ with $g\cap\sigma=\ell \cap \sigma$. 

		Now we consider the cases with $\ell\cap \sigma\neq g\cap \sigma$. First, we consider the line-plane pairs $(g,\tau)$ who satisfy $g\cap \pi\neq \emptyset$.
		The line $g$ is determined by one of the $q^2+q$ points $\neq \ell\cap \sigma$ on $\sigma$ and one of the $q^2+q$ points $\neq \ell\cap \sigma$ on $\pi$.
 Every line $g$ is incident with $\lambda$ planes $\tau$. This yields $(q^2+q)^2\lambda$ coflags.

		 Now we look at the chambers  $(g,\tau)$ who satisfy $\ell\cap \tau\neq \emptyset$. We can choose one of $q^2+q$ points $\neq l\cap \sigma$ for $\tau\cap \sigma$ on $\sigma$ and one of $q$ points on $\ell$ that is different from $\ell\cap \sigma$, for $\ell\cap \tau$.
		  These two points span a line that has to be incident with $\tau$, therefore $\tau$ is one of the $\lambda$ planes that is incident with this line. The line $g$ has to be one of the $q+1$ lines on $\tau$ that is incident with $\tau\cap \sigma$. This yields $(q^2+q)q(q+1)\lambda$ coflags.	

In conclusion $(\#a)\leq(q^3+q^2)\lambda+(q^2+q)^2\lambda+(q^2+q)q(q+1)\lambda= \lambda(2q^4+5q^3+3q^2)$.

	\item[b)] We find an upper bound for the number $(\#b)$ of chambers $(Q,g,\tau,R)$ in $Y_1$, whose point $Q$ is in $S$.

Since $\pi\cap \sigma$ is a point and $\pi\subseteq S$, the subset $S\cap \sigma$ is a line.

If $(Q,g,\tau,R)$ is a chamber in $Y_1$ with $Q\in S$, then $Q$ has to be on the line $S\cap g$.
Therefore $g$ has to meet the line $S\cap \sigma$. There are $q+1$ points on the line $S\cap\sigma$. Each point is incident with $q^3+q^2$ lines not in $\sigma$ and each line is incident with $\lambda$ planes. Since the coflags of $Y_1$ have weight one, we have $(\#b)\leq (q+1)(q^3+q^2)\lambda$.

	\item[c)] We want to find an upper bound for the number $(\#c)$ of chambers $(Q,g,\tau,R)$ in $Y_1$, whose solid contains $P$.

 Recall that we only consider chambers $(Q,g,\tau,R)$ whose solids do not contain $\sigma$, since otherwise $\tau\cap \sigma$ is a line. Therefore $R\cap \sigma$ has to be one of the $q+1$ lines in $\sigma$ incident with $P$. 

	If all chambers $(Q,g,\tau,R)$ of $Y_1$ that satisfy $P\in R$, also contain a  line-plane pair $(g,\tau)$ that is not fopposite to $(\ell,\pi)$, we have $(\#c)\leq (\#a)$.

Now assume that there is a chamber $(P_1,\ell_1,\pi_1,S_1)$ in $Y_1$ with the following properties: $P\in S_1$ and $(\ell,\pi)$ and $(\ell_1,\pi_1)$ are f-opposite.
We find an upper bound for the number of chambers $(Q,g,\tau, R)$ in $Y_1$ that satsify $P\in R$ and that are non-opposite to $(P_1,\ell_1,\pi_1,S_1)$ and $(P,\ell,\pi,S)$. This will give an upper bound for $(\#c)$.

		\begin{itemize}
		\item The number of chambers $(Q,g,\tau, R)$ who satisfy $P\in R$ and $(g,\tau)$ is non-f-opposite to $(\ell_1,\pi_1)$ is at most $\lambda(2q^4+5q^3+3q^2)$. This follows directly from a) applied to $(\ell_1,\pi_1)$. 
		\item The number of chambers $(Q,g,\tau,R)$ who satisfy $P\in R$ and $Q\in S_1$ is at most $(q+1)(q^3+q^2)\lambda$. This follows directly from b) applied to $(P_1,\ell_1,\pi_1,S_1)$. 
		\item Now we consider the chambers $(Q,g,\tau,R)$ who satisfy $P\in R$ and $P_1\in R$. Since $(\ell,\pi)$ and $(\ell_1,\pi_1)$ are f-opposite we have $P\neq P_1$. The set $R\cap \sigma$ has to be a line, since $R$ cannot contain $\sigma$. 
			There are $q+1$ points $Q$ on $PP_1$.
						Every point is incident with $q^3+q^2$ lines $g$ not in $\sigma$. Every line is incident with $\lambda$ planes $\tau$. 
 Therefore, we have at most $(q+1)(q^3+q^2+q+1)\lambda$ coflags
		\end{itemize}
		This yields 
		\begin{align*}
			\begin{split}
				(\#c)\leq & (2q^4+5q^3+3q^2)\lambda\\
				&+(q+1)(q^3+q^2)\lambda\\
				&+ (q+1)(q^3+q^2)\lambda.
			\end{split}
		\end{align*}
\end{itemize}

In conclusion we have
\begin{align*}
	\begin{split}
		|Y_1|\leq(\#a)+(\#b)+(\#c) &\leq 2( 2q^4+5q^3+3q^2)\lambda \\
		&+ 3(q+1)(q^3+q^2)\lambda\\
		&=   (7 q^4 + 16 q^3 + 9 q^2)\lambda.
	\end{split}
\end{align*}
Considering $|Y|=|Y_3|+|Y_2|+|Y_1|$, we obtain
 \begin{align*}
	\begin{split}
		|Y| &\leq  2(q^2+q+1)\lambda \\
		&+  (q^2+q+1)(q+1)^2\lambda \\
		&+   (7 q^4 + 16 q^3 + 9 q^2) \lambda \\
		&=  (8 q^4 + 19 q^3 + 15 q^2 + 5 q + 3)\lambda.
	\end{split}
\end{align*}
\end{proof}

\begin{lemma} \label{L: points in a solid preparation CASE A} 
Let $M$ be an independent set of $\Gamma_4$. Let $\sigma$ be a plane and let $Y$ be the set of all chambers $(P,\ell,\pi,S)$ of $M$ for which $(\ell,\pi)$ has $M$-weight at most two and for which $P\in \sigma$. Furthermore, assume that we are in Case A. We have $|Y|\leq y_A:=6q^5+12q^4+8q^3+4q^2$.
\end{lemma}

\begin{proof}
First consider Case A1.
	There are $q^2+q+1$ points on $\sigma$. For every point of $\sigma$ the number of lines $\ell$ that meet $\sigma$ is at most $q^3+q^2+q+1$. 
	Let $\ell$ be a line of a chamber $(P,\ell,\pi,S)$ in $Y$. Then Remark \ref{R: raez} implies that $\pi$ is uniquely determined by $\ell$ in Case A1.
	Since every $M$-coflag has weight at most two, we get $|Y|\leq 2(q^2+q+1)(q^3+q^2+q+1)$.
	
	Now consider Case A2. In this case all coflags $(\ell,\pi)$ of $Y$ contain lines $\ell$ that have to be in a fixed solid $S_0$ and planes $\pi$ that satisfy $\pi\cap S_0=\ell$. 	This implies that every line $\ell$ is incident with at most $q^2$ planes $\pi$. If $\sigma\cap S_0$ is just a line, there are $(q+1)(q^2+q)+1$ lines in $S_0$ that meet $\sigma\cap S_0$ and therefore $|Y|\leq 2((q+1)(q^2+q)+1)q^2$. Now consider the case $\sigma \subseteq S_0$.
	
	Let $Y=Y_1\dot{\cup} Y_2$, where $Y_2$ contains all chambers with a line in $\sigma$. Then $|Y_2|\leq 2(q^2+q+1)q^2$.
	
	Similar to the proof of the previous lemma we get the following:
	Since the point of every chamber of $Y_1$ is in $\sigma$, but the lines of the chambers of $Y_1$ are not contained in $\sigma$, Proposition \ref{P: weight of lines} yields we can count the $Y_1$-coflags in order to count the chambers in $Y_1$.

	Let $(P_1,\ell_1,\pi_1,S_1)$ be a chambers in $Y_1$. First, we find an upper bound for the number $n$ of chambers $(Q,g,\tau,R)$ in $Y_1$ whose coflag is non-f-opposite to $(\ell_1,\pi_1)$. \\
	 If $\ell_1$ and $\tau$ meet, we have that $\ell_1$ and $g$ meet, since $\tau\cap S_0=g$.  If $\pi_1$ and $g$ meet, we also have that $\ell_1$ and $g$ meet, since $\pi\cap S_0=\ell$.
	 There are $q+1$ points on $\ell_1$. Every point is incident with $q^2+q+1$ lines $g$ in $S_0$. The plane $\tau$ is one of the $q^2$ planes incident with $g$ that is not in $S_0$. In conclusion $n\leq (q+1)(q^2+q+1)q^2$.
	 	
	Let $(Q_1,g_1,\tau_1,R_1)$ be a chamber in $Y_1$ with $(g_1,\tau_1)$ f-opposite to $(\ell_1,\pi_1)$ and without loss of generality we assume that the line $R_1\cap \sigma$ contains $P_1$. Note that $S_1\cap \sigma$ is also a line. Let $(Q,g,\tau,R)$ be another chamber of $Y_1$ and assume: $(g,\tau)$ is f-opposite to $(\ell_1,\pi_1)$ and $(g_1,\tau_1)$, and $Q\notin R_1\cap \sigma,  S_1\cap \sigma$. In this case $R$ has to contain $Q_1$ and $P_1$ which implies $Q\in R_1\cap \sigma$, contradiction.
	Therefore all chambers in $Y_1$ have to be non-f-opposite to $(\ell_1,\pi_1)$, or non-f-opposite to $(g_1,\tau_1)$, or they have to have their point on one of the $2q+1$ points of $h:=(R_1\cap \sigma) \cup (S_1\cap \sigma)$. There are at most $(2q+1)(q^2+q)q^2$ chambers in $Y_1$ with a point on $h$.\\ This implies $|Y_1|\leq 2n+ 2(2q+1)(q^2+q)q^2$.\\ Using the bounds for $|Y_1|, |Y_2|$ and $Y=Y_1\dot{\cup} Y_2$ yields the desired bound.
\end{proof}

	\begin{lemma} \label{L: points in a solid} 
	Let $M$ be an independent set of $\Gamma_4$. Let $S$ be a solid and let $X$ be the set of all chambers $(P,\ell,\pi,R)$ of $M$ for which $(\ell,\pi)$ has $M$-weight at most two and for which $P\in S$. \\
	In Case A2 we also assume $S\neq S_0$.\\
	 In Case A we have $|X|\leq x_A:=30 q^5 + 65 q^4 + 58 q^3 + 44 q^2 + 17 q + 8$.\\
	 In Case B we have $|X|\leq x_B:= 184 q^6 + 583 q^5 + 959 q^4 + 959 q^3 + 676 q^2 + 300 q + 116$.
\end{lemma}

\begin{proof}
The points of the chambers in $X$ will be called $X$-points. 

In Case A2 we have that $S\cap S_0$ is a plane, since $S\neq S_0$. Therefore, the result follows directly from Lemma \ref{L: points in a solid preparation CASE A}. Now consider Case A1 and B simultaneously. 	If we are in Case A1 let $i:=i_A$, and if we are in Case B let $i:=i _B$ for $i\in \{z,y\}$.
	If we are in Case A1 let $\lambda:=1$, if we are in Case B let $\lambda:=q^2+q+1$.

Let $X_3$ be the set of chambers in $X$ whose line is contained in $S$. Let $X_3^c$ be the set of chambers in $X$, whose line is not contained in $S$. 
Clearly $X= X_3^c \cup X_3 $. Since we want to obtain bounds for $|X_3|$ and $|X_3^c|$ we assume without loss of generality that $X_3$ and $X_3^c$ are not empty.\\
There are $(q^2+1)(q^2+q+1)$ lines in $S$ and each line is incident with $\lambda$ planes, therefore the number of distinct $X_3$-coflags is at most $(q^2+1)(q^2+q+1)\lambda$ and hence $|X_3|\leq 2(q^2+1)(q^2+q+1)\lambda$.\\

Now we find an upper bound for $|X_3^c|$. 
Since the point of every chamber of $X_3^c$ is in $S$, but the lines of the chambers of $X_3^c$ are not contained in $S$,  we can count the $X_3^c$-coflags, instead of the chambers in $X_3^c$.\\
If all $X_3^c$-points lie in a plane $\sigma \subseteq S$, Lemma \ref{L: points in a solid preparation} and Lemma \ref{L: points in a solid preparation CASE A} yield $|X_3^c|\leq y$ and the assertion follows. Therefore, we may assume that the $X_3^c$-points span $S$.\\

	Let $P_1,P_2,P_3,P_4$ be four $X_3^c$-points that span the solid $S$. Let $(P_i,\ell_i,\pi_i,S_i)$ be a chamber in $X_3^c$ for $i=1,\ldots,4$. \\
	The number of chambers $(Q,g,\tau,R)$ in $X_3^c$ with a point in the span $\erz{ P_1,P_2,P_3 }$ is, according to Lemma \ref{L: points in a solid preparation}, at most $y$.\\
	For $i=1,2,3$ the line $\ell_i$ is not in $S$, therefore $S_i\cap S$ is a plane. Hence the number of chambers in $X_3^c$, whose point is in $S_1\cup S_2 \cup S_3$, is at most $3\cdot y$.\\

Let $X_1$ be the set of chambers in $X_3^c$, whose points are not in $\erz{ P_1,P_2,P_3 }$, $S_1$, $S_2$, or $S_3$. We now obtain an upper bound for $|X_1|$, therefore we assume without loss of generality that $X_1$ is not empty. Since $|X_3^c|\leq 4y+|X_1|$, this provides an upper bound for $|X_3^c|$.\\

Let $(Q,g,\tau, R)$ be a chamber in $X_1$. Since $g$ is not in $S$ and $Q$ is not in $\erz{P_1,P_2,P_3}$, or $S_1\cup S_2\cup S_3$ , the solid $R$ cannot contain $P_1$, $P_2$ and $P_3$.  Therefore the $X_1$-coflag $(g,\tau)$ has to be non-f-opposite to one of the coflags $(\ell_i,\pi_i)$ for $i=1,2,3$. \\
We find an upper bound for the number of chambers in $X_1$, whose line-plane pair is non-f-opposite to $(\ell_1,\pi_1)$. Multiplying this bound by three, yields the desired bound for $|X_1|$.\\

First consider the Case A1. If $(g,\tau)$ is non-f-opposite to $(\ell_1,\pi_1)$, we have that $g$ meets $\pi_1$ or that $\ell_1$ meets $\tau$. The number of lines $g$ that meet $\pi_1$ is at most $(q^2+q+1)(q^3+q^2+q+1)$. Finally, $\tau$ is determined by $g$ and $P_0$. \\
Now consider one of the $q+1$ points on $\ell_1$. The line that is spanned by this point and $P_0$ has to be incident with $\tau$, that is $\tau$ is one of the $q^2+q+1$ planes incident with this line. Finally, $g$ is one of the $q^2+q+1$ lines on $\tau$. In conclusion we have 
\begin{align*}
\begin{split}
				|X_1|&\leq \  3\cdot \Bigl(  
							(q^2+q+1)(q^3+q^2+q+1)+(q+1)(q^2+q+1)^2 \Bigr) \\
							&=6 q^5 + 15 q^4 + 24 q^3 + 24 q^2 + 15 q + 6.
	\end{split}
	\end{align*}

Now consider Case B.
To obtain our bound for the number of chambers $(Q,g,\tau,R)$ in $X_1$, whose line-plane pair $(g,\tau)$ is non-f-opposite to $(\ell_1,\pi_1)$, we study the three different possible relations between $(Q,g,\tau,R)$ and $(P_4,\ell_4,\pi_4,S_4)$ and obtain an upper bound for each case. Recall that all the lines $\ell_1,\ldots,\ell_4,g$ are not in $S$ and that $\pi_i\cap S$ and $\tau\cap S$ are lines. Also recall that all coflags in $X_1$ have weight one.
	\begin{itemize}

		\item[a)] Assume $Q\in S_4$ and $(g,\tau)$ is non-f-opposite to $(\ell_1,\pi_1)$.
		
		 Since $S_4\cap S$ is a plane, Lemma \ref{L: points in a solid preparation} and Lemma \ref{L: points in a solid preparation CASE A} yield that at most $y$ chambers $(Q,g,\tau,R)\in X_1$ satisfy this condition.

		\item[b)] Assume $P_4\in R$ and $(g,\tau)$ is non-f-opposite to $(\ell_1,\pi_1)$.
		
				The number of line plane-pairs $(g,\tau)$, with $P_4\in g$, is at most $(q^3+q^2+q+1)\lambda$. From now on we assume $P_4\notin g$.\\
				Let $(Q_1,g_1,\tau_1,R_1)$ be a chamber with $P_4\notin g_1$, $P_4\in R_1$. There are at most $2(q^2+q+1)$ such chambers with $g_1=\ell_1$, or $\tau_1=\pi_1$, from now on we assume $g_1\neq \ell_1$ and $\tau_1\neq \pi_1$. All chambers $(Q,g,\tau,R)$ with $P_4\in R\setminus g$ and $(g,\tau)$ non-f-opposite to $(\ell_1,\pi_1)$, have to satisfy one of the three following conditions: $Q\in R_1$, or $(g,\tau)$ is non-f-opposite to $(g_1,\tau_1)$, or $Q_1\in R$.
				
				The number of chambers that satisfy $Q\in R_1$ is at most $y$, since $R_1\cap S$ is a plane.
				
				Lemma \ref{L: z 2 line-plane flags opposite} and Lemma \ref{L: z 2 line-plane flags opposite, Case A} yield that the number of chambers that are non-opposite to $(g_1,\tau_1)$ and $(\ell_1,\pi_1)$ is at most $z$, since $g_1\neq \ell_1$ and $\pi_1\neq \tau_1$.
				
				Finally, we want to find an upper bound for the number of chambers $(Q,g,\tau,R)$ that satisfy $P_4\in R$ and $Q_1\in R$. Since we assumed $P_4\notin g_1$, the points $P_4$ and $Q_1$ span a line. In Lemma \ref{L: points in a solid preparation}, we showed that the number of chambers, whose point is in a given plane, is at most $y$. Dually the number of chambers whose solid contains a given line, is at most $y$

		\item[c)] Assume $(g,\tau)$ is non-f-opposite to $(\ell_4,\pi_4)$ and $(g,\tau)$ is non-f-opposite to $(\ell_1,\pi_1)$.

Since $\ell_1$ and $\ell_4$ are not in $S$ and $P_1\neq P_4$, we have $\ell_1\neq \ell_4$.
If $\pi_4\neq \pi_1$, Lemma \ref{L: z 2 line-plane flags opposite} implies that at most $z$ line-plane flags $(g,\tau)$ satisfy the assumed condition. Therefore, let us assume $\pi_4=\pi_1$.

The number of line-plane pairs $(g,\tau)$ with $\tau=\pi_1$ is at most $q^2+q+1$, from now on consider only coflags $(g,\tau)$ with $\tau\neq \pi_1$.\\
			Since $(Q,g,\tau,R)\in X_3^c$, we have $Q\notin S_2\cup S_3$.
		 
		 If $(Q,g,\tau,R)$ satisfies $P_2,P_3\in R$, the line $P_2P_3$ meets $\tau$ in a point. Note that $\pi_1=\pi_4$ implies $\pi_1\cap S=P_1P_4$, hence $P_2P_3\cap \pi_1=\emptyset$.

First, we obtain an upper bound for the number of line-plane pairs $(g,\tau)$ with $g\cap \pi_1\neq \emptyset$ and $P_2P_3\cap \tau \neq \emptyset$: We can choose one of the $q^2+q+1$ points on $\pi_1$ for $g\cap \pi_1$ and one of the $q+1$ points on $P_2P_3$ for $P_2P_3\cap \tau$. The line spanned by these two points is incident with $\lambda$ planes $\tau$ and $g$ has to be one of the $q+1$ lines on $\tau$ that is incident with $g\cap \pi_1$. Hence the desired bound is $(q^2+q+1)\lambda(q+1)^2$.

Now we obtain an upper bound for the number of line-plane pairs $(g,\tau)$, with $\tau \cap \ell_1, \tau\cap \ell_4 \neq \emptyset$ and $P_2P_3\cap \tau \neq \emptyset$:
Since $\ell_1\cap S=P_1$ and $P_2P_3\subseteq S$, we have $\ell_1\cap P_2P_3=\emptyset$.
 If $\tau\cap \ell_1\neq\tau\cap \ell_4$, we choose one of the $q+1$ points on $\ell_1$ 
 one of the $q+1$ points on $P_2P_3$. These two points span a line that is incident with  $\tau$ and $g$ is one of the $q^2+q+1$ lines on $\tau$.\\ 
 If $\tau\cap \ell_1=\tau\cap \ell_4$, we choose any of the $q+1$ points on $P_2P_3$ and the line spanned by this point and $\tau\cap \ell_1$ is incident with $\tau$. Hence $\tau$ is one of the $\lambda$ planes incident with this line and $g$ is one of the $q^2+q+1$ lines on $\tau$. In conclusion our bound is $(q+1)^2 \lambda (q^2+q+1)+(q+1)(q^2+q+1) \lambda$.

Thus far we have only considered line-plane pairs $(g,\tau)$ that are non-f-opposite to $(\ell_1,\pi_1)$ and $(\ell_4,\pi_4)$ and that occur in a chamber $(Q,g,\tau , R)$ that satisfies one of the following conditions: $Q\in S_2$, $Q\in S_3$, $P_2$ and $P_3\in R$. In this case \\
Now consider a line plane pair $(\tilde{g},\tilde{\tau})$ that is non-f-opposite to $(\ell_1,\pi_1)$ and $(\ell_4,\pi_4)$ and that occurs in a chamber $(Q,\tilde{g},\tilde{\tau} , R)$ that does not satisfy any of the following conditions: $Q\in S_2$, $Q\in S_3$, $P_2$ and $P_3\in R$. Then $(\tilde{g},\tilde{\tau})$ has to be non-f-opposite to $(\ell_2,\pi_2)$, or $(\ell_3,\pi_3)$. 
Since $\pi_1=\pi_4$, we have $\pi_4\neq \pi_2,\pi_3$. The points $P_2,P_3,P_4$ are in $S$, but the lines $\ell_2,\ell_3,\ell_4$ are not in $S$, therefore the lines are pairwise distinct. Lemma \ref{L: z 2 line-plane flags opposite} and Lemma \ref{L: z 2 line-plane flags opposite, Case A} imply that the number of line-plane pairs that are non-opposite to $(\ell_4,\pi_4)$ and non-opposite to $(\ell_2,\pi_2)$  or $(\ell_3,\pi_3)$ is at most $2z$.

\end{itemize}
In conclusion we obtain in Case B:
\begin{align*}
	\begin{split}
			|X_1|\leq \ & 3\cdot \Bigl(  y \\
							&+(q^3+q^2+q+1)\lambda+2(q^2+q+1)+y+z+y \\
							&+(q^2+q+1)
							+(q^2+q+1)\lambda(q+1)^2
							+(q+1)^2 \lambda (q^2+q+1)+(q+1)(q^2+q+1)\lambda
							+2z \Bigr)  \\
							&=150 q^6 + 471 q^5 + 783 q^4 + 795 q^3 + 576 q^2 + 264 q + 102
	\end{split}
\end{align*}
Using $|X| \leq |X_3|+4y+|X_1|$ and the bounds obtained for $|X_3|$ and $|X_1|$, we get the desired result.
\end{proof}

\begin{prop} \label{P: 1 and 2}
Let $M$ be an independent set of $\Gamma_4$. And let $W$ be the set of all all chambers in $M$ that contain a coflag of weight at most two.\\
In Case A, we have $|W|\leq 122 q^5 + 272 q^4 + 243 q^3 + 184 q^2 + 71 q + 35$.\\
In Case B, we have $|W|\leq 744 q^6 + 2354 q^5 + 3874 q^4 + 3876 q^3 + 2736 q^2 + 1216 q + 470$
\end{prop}

\begin{proof}
If we are in Case A, let $x:=x_A$. If we are in Case B let $x:=x _B$.

	\underline{Case 1:} No two $M$-coflags of weight at most two are f-opposite.\\
 In this case Theorem \ref{T: B and B q^5} yields that the number of $M$-coflags with weight at most two, is at most $q^5+3q^4+4q^3+4q^2+2q+1$. This implies $|W|\leq 2(q^5+3q^4+4q^3+4q^2+2q+1)$.\\
	
	\underline{Case 2:} There are two $M$-coflags $(\ell_1,\pi_1)$ and $(\ell_2,\pi_2)$ of weight at most two that are f-opposite.\\
	Let $(P_1,\ell_1,\pi_1,S_1)$ and $(P_2,\ell_2,\pi_2,S_2)$ be chambers of $M$. All chambers $(Q,g,\tau, R)$ in $M$, with $M$-coflags of weight at most two, must satisfy one of the following conditions:
	\begin{itemize}
		\item[(i)] $(g,\tau)$ is non-f-opposite to $(\ell_1,\pi_1)$ and non-f-opposite to $(\ell_2,\pi_2)$. 
		\item[(ii)] $Q\in S_1$.  
		\item[(iii)] $P_1\in R$.
		\item[(iv)] $Q\in S_2$. 
		\item[(v)] $P_2\in R$.
	\end{itemize}
	Since $(\ell_1,\pi_1)$ and $(\ell_2,\pi_2)$ are f-opposite, we have $\ell_1\neq \ell_2$ and $\pi_1\neq \pi_2$. Therefore,  Lemma \ref{L: z 2 line-plane flags opposite} and \ref{L: z 2 line-plane flags opposite, Case A} yield that the number of chambers that satisfy (i) is at most $z$.\\
	Consider Case A2. Since $(\ell_1,\pi_1)$ has weight at most two, the plane $\pi_1$ is not contained in $S_0$. This implies that in Case A2 any solid  $S_1$ incident with $\pi_1$ is distinct from $S_0$.\\
	 Therefore, the number of chambers that satisfy any of the cases (ii)-(v) is bound by Lemma \ref{L: points in a solid}, or the dual statement respectively. This yields that the number of $M$-chambers with $M$-coflags of weight one or two is $\leq 4x+z$.
\end{proof}

\begin{cor}
Let $M$ be a maximal independent set of $\Gamma_4$ that contains coflags of weight $<(q+1)^2$. \\
In Case A we have $|M|\leq q^7 + 4 q^6 + 133 q^5 + 290 q^4 + 261 q^3 + 195 q^2 + 75 q + 36$. \\ 
In Case B we have $|M|\leq 824 q^6 + 2551 q^5 + 4362 q^4 + 4605 q^3 + 3589 q^2 + 1783 q + 710$.
\end{cor}
\begin{proof}
Proposition \ref{P: weight of lines} implies that a coflag can have weight $0$, $1$, $2$, $q+1$, $2q+1$, or $(q+1)^2$.
 We utilize the Propositions \ref{P: upper bound for (q+1)^2-lines}, \ref{P: q+1}, \ref{P: q+1 CASE A}, and \ref{P: 1 and 2} an obtain the following bounds.\\
 In Case A we have 
 \begin{align*}
 	\begin{split}
		|M|& \leq ((q^3+q^2+q+1)(q^2+q+1)+q^2(q+1))(q+1)^2 \\
		&+ 2\cdot (q^4+q^3+q^2)(q+1)\\
		&+ 122 q^5 + 272 q^4 + 243 q^3 + 184 q^2 + 71 q + 35 \\
		&=q^7 + 4 q^6 + 133 q^5 + 290 q^4 + 261 q^3 + 195 q^2 + 75 q + 36.
	\end{split}
 \end{align*}
 In Case B we have 
 \begin{align*}
 \begin{split}
		|M|& \leq (4q^4+9q^3+4q^2+q+1)(q+1)^2 \\
		&+ 2\cdot(q^5+3q^4+4q^3+4q^2+2q+1) (q+1)\\
		&+ 744 q^6 + 2354 q^5 + 3874 q^4 + 3876 q^3 + 2736 q^2 + 1216 q + 470 \\
		&=750 q^6 + 2379 q^5 + 3914 q^4 + 3910 q^3 + 2755 q^2 + 1225 q + 473.
	\end{split}
 \end{align*}
\end{proof}

This completes the proof of Theorem \ref{T: main theorem}.



\subsection*{Acknowledgements}
  
The author would like to thank Klaus Metsch for bringing this problem to his attention and for providing a lot of helpful discussions.

\bibliographystyle{plainurl}
\bibliography{Maximal_independent_sets_of_chambers_of_PG4q_Heering_arxiv}

\end{document}